%% file: CDC_RSTracking.tex
\documentclass[12pt]{article}
\usepackage{latexsym}
\usepackage{amsmath}
\usepackage{amssymb}
\usepackage{amsthm}
\usepackage{stmaryrd}
\usepackage{epsfig}
\usepackage{graphicx}
\usepackage{enumerate}
\usepackage{subfigure}
\usepackage{url}

\theoremstyle{plain}
\newtheorem{thm}{Theorem}[section]

\newtheorem{proposition}[thm]{Proposition}

\theoremstyle{definition}

\theoremstyle{remark}
\newtheorem*{remark}{Remark}

\newcommand{\be}{\begin{enumerate}}
\newcommand{\ee}{\end{enumerate}}

\newcommand{\Tr}{\text{Tr}}

\advance \topmargin by -\headheight
\advance \topmargin by -\headsep
\evensidemargin \oddsidemargin
\usepackage[pdftex,
    pdftitle={Performance Evaluation of a Multi-Agent Risk-Sensitive Tracking System},
    pdfauthor={Jerome Le Ny},
    colorlinks,linkcolor={blue},citecolor={blue},urlcolor={red},
]
{hyperref}

\title{\LARGE \bf
Performance Evaluation of a Multi-Agent Risk-Sensitive Tracking System
}

\author{ \parbox{3 in}{\centering Jerome Le Ny and Eric Feron}
         \thanks{This work was supported by Air Force - DARPA - MURI award 009628-001-03-132 and Navy ONR award N00014-03-1-0171. J. Le Ny is with the Laboratory for Information and Decision Systems, Massachusetts Institute of Technology, Cambridge, MA 02139-4307, USA. E. Feron is with the School of Aerospace Engineering, Georgia Tech, Atlanta, GA 30332, USA. {\tt\small jleny@mit.edu, eric.feron@aerospace.gatech.edu} }\\
}

\begin{document}

\maketitle
\thispagestyle{empty}
\pagestyle{empty}

\begin{abstract}

In this paper, we consider a simple linear exponential quadratic Gaussian (LEQG) tracking problem for a multi-agent system. We study the dynamical behaviors of the group as we vary the risk-sensitivity parameter, comparing in particular the risk averse case to the LQG case. Then we consider the evolution of the performance per agent as the number of agents in the system increases. We provide some analytical as well as simulation results. In general, more agents are beneficial only if noisy agent dynamics and/or imperfect measurements are considered. The critical value of the risk sensitivity parameter above which the cost becomes infinite increases with the number of agents. In other words, for a fixed positive value of this parameter, there is a minimum number of agents above which the cost remains finite.

\end{abstract}

\section{INTRODUCTION}

With the rapidly growing interest in sensor networks and distributed control systems, a large body of litterature in recent years has focused on new and old engineering problems that these networks pose. Examples include consensus problems in various forms (see \cite{Blondel05_consensus} and the references therein), control design with pre-specified information passing structure and control under communication constraints (e.g. \cite{Tatikonda00_controlCom, langbort04_distributedControl}),  designing communication schemes in ad-hoc networks \cite{Akyildiz02_survey}, etc.

In some cases, such as in signal processing and control for distributed parameter systems (e.g. \cite{Willsky02_MMR, Bahmieh02_distributed}), the distributed architecture of the system is imposed by the task in mind. But in a large number of situations, especially involving mobile robotic networks, the designer has a significant freedom in choosing the degree of distributivity of the system, the number and type of agents to use, etc. Yet relatively little work has focused on understanding quantitatively the benefits of multi-agent systems over centralized systems for various tasks, or the impact of the number of agents on the performance of the system, when this choice is available.

For a given task at hand, it is important to study these problems, because they influence most of the system design. For example, scalability is a relative notion. If we find that no significant improvement, or worse a decrease in performance is to be expected if we employ thousands of agents instead of tens or hundreds, it is probably better to design algorithms with good performance that are not necessarily the most decentralized. Moreover, it is clear than more agents is not always beneficial. For example, the capacity of a wireless sensor network is expected to decrease as the number of nodes increases \cite{Gupta00_capacity}. A similar issue arises in mobile robotic networks where conflict free vehicle routing is necessary and congestion increases with the size of the network \cite{savchenko05_complexity}. 

In some recent work researchers have obtained asymptotic performance scalings for certain multi-agent systems. Let us mention \cite{enright05_multiUAV, waisanen06_delay} on various dynamic routing problems. One can also find in the computer science litterature a large number of papers focusing on the minimum number of agents necessary to perform certain tasks: for instance, the minimum number of pursuers to catch an evader \cite{isler06_pursuitEvasion}, or the minimum number of guards for an art gallery \cite{Shermer92_artGallery}.

In this work, we study a simple multi-agent tracking problem using linear quadratic control tools, which are well suited for performance evaluation. We focus on the evolution of the performance per agent, as the number of agents increases. Numerous papers describe multi-agent architectures for tracking and estimation, the simultaneous localization and mapping problem, etc., focusing on the various subroutine designs described earlier (e.g. \cite{Martinez06_tracking}). Because of the link between risk-sensitive control and differential games, our application is also related to some recent work on multi-agent pursuit-evasion games \cite{Hespanha00_pursuit}. The model we consider is simpler, but our focus is on obtaining insight regarding the performance asymptotics of the system, which these references do not usually discuss. 

The rest of the paper is organized as follows. Section \ref{section: model} describes our model and section \ref{section: solution} the general form of its solution. We find in particular that a standard LQG formulation might not provide a satisfying cooperative control law for the group of agents, as the individual controllers decouple, essentially because of the certainty equivalence principle. We can obtain coupled control laws however if we consider risk-sensitive agents. Intuitively, we expect the set of risk-sensitive parameters where tracking is possible, or the robustness to model disturbances of our system, to increase with the number of agents. This is related to the earlier observation that a minimum number of agents can be necessary in certain pursuit-evasion games. Section \ref{section: performance} provides some analysis of this aspect. Increasing the number of agents becomes critical in our scenario only in the case of noisy dynamics and, of course, noisy measurements since more agents obtain a better state estimator. We give elements of analytical performance analysis as well as simulation results in that section as well.

\input{model}

\vspace{0.1cm}

\input{risk}

\input{performance}

\section{CONCLUSION}
In this work we have considered a basic tracking task to be performed cooperatively by a team of mobile sensors. We have given some elements of analysis concerning the influence of the size of the group on the individual performance. Intuitively, there seems to be much gain to be obtained by multi-agent systems in terms of robustness and in a risk-sensitive context, and we believe that more work in this direction is needed. In any case, it is clear that there is a need for a better understanding of the beneficial role of distributed architectures, when this aspect is under the control of the system designer. 



\input{thebiblio}
\end{document}

%% file: model.tex
\section{Basic Model} \label{section: model}

We will consider the following tracking problem. 
There is an evader moving randomly in $\mathbb{R}^d$, subject to the linear dynamics
\begin{equation}	\label{eq: evader's dynamics}
\dot x_e = A x_e + G w_e,
\end{equation}
where $w_e$ is a standard $d$-dimensional white Gaussian noise, and we let $W=GG'$. 

We have $n$ identical pursuers, which are also described by linear systems, with the same $A$ matrix for simplicity: 
\[
\dot x_{p,i} = A x_{p,i} + B u_{p,i} + \sqrt \epsilon F w_{p,i}, \quad i=1,\ldots,n,
\]
where $w_{p,i}$ are standard $d$-dimensional white Gaussian noises, independent between different agents and from $w_e$. $\epsilon$ is a parameter that will tend to $0$ in parts of the analysis later on. Let $Z=F F'$. 
The initial positions of the agents are also independent d-dimensional Gaussian random vectors.

Each pursuer incurs a running cost 
\begin{equation} \label{eq: state cost}
\frac{1}{2} (x_{p,i}-x_e)' Q (x_{p,i}-x_e)
\end{equation}
which is quadratic in the difference between its state and the state of the evader. There is also a running control cost $\frac{1}{2} u_{p,i}' \, R \, u_{p,i}$ for each agent, where $R$ is a positive definite matrix. Alternatively, the $n$ mobile agents are tracking the center of an ellispoid moving randomly in $\mathbb{R}^d$ according to (\ref{eq: evader's dynamics}), whose shape is known and given by a positive semi-definite matrix $Q$. A motivation for this model is cooperative soaring and tracking ascending currents for a team of Unmanned Aerial Vehicles.

Each agent has access to a relative measurement:
\[
y_{p,i}= C (x_{p,i}-x_e) + H v_{p,i}, \; i=1,\ldots,n,
\]
with $v_{p,i}$ a standard d-dimensional white noise. We let $V=HH'$ and assume $V$ to be positive definite. The measurement noise processes of the different agents are independent, and also independant of the various noises in the dynamics. For example, the agents could obtain noisy measurements of the gradient of the quadratic cost function ($\ref{eq: state cost}$), in which case $C=Q$. We will also consider the \emph{perfect measurement case}, where $y=x$, subject to no measurement noise.

We define $x_{i}=x_{p,i}-x_e$, and the aggregate vectors $x=[x'_1,\ldots,x'_n]'$, $y=[y'_{p,1},\ldots,y'_{p,n}]'$, $u=[u'_{p,1},\ldots,u'_{p,n}]'$, $w_p=[w'_{p,1},\ldots,w'_{p,n}]'$, $v=[v'_{p,1},\ldots,v'_{p,n}]'$. We will also use the Kronecker product of matrices, which we recall. If $A=[a_{ij}]$ and $B$ are matrices, then we have by definition:
\[
A \otimes B = 
\left[ \begin{array}{ccc}
a_{11}B & \hdots & a_{1n}B \\
\vdots & \ddots & \vdots \\
a_{m1}B & \hdots & a_{mn}B
\end{array}
\right] \; ,
\] 
and we will use the property:
\[
(A \otimes B)(C \otimes D)=AC \otimes BD.
\]
The eigenvalues of $A \otimes B$ are $\{\lambda_i \mu_j\}$ and the corresponding eigenvectors $\{x_i \otimes y_j\}$, where $\{\lambda_i\}$, $\{x_i\}$, $\{\mu_j\}$, $\{y_j\}$ are the eigenvalues and eigenvectors of $A$ and $B$ respectively. Hence if $A$ and $B$ are positive definite matrices, so is $A \otimes B$. Denote by $I_n$ the $n \times n$ identity matrix, by $\mathbf{1}_n$ the column vector of ones of size $n$, and by $E_n=\mathbf{1}_n \mathbf{1}'_n$ the $n \times n$ matrix of ones. Let 
\begin{align*}
&A_n=I_n \otimes A, \; B_n=I_n \otimes B, \; F_n=I_n \otimes F, \; Z_n=I_n \otimes Z, \\
&C_n=I_n \otimes C, \; H_n=I_n \otimes H, \; V_n= I_n \otimes V,
\end{align*}
be the block diagonal matrices describing the dynamics of the group of $n$ pursuers. Finally, let $G_n=\mathbf{1}_n \otimes G$ and the corresponding spectral density matrix be 
\[
W_n=E_n \otimes W=G_n G'_n= \left[ \begin{array}{ccc} W & \cdots & W \\ \vdots & & \vdots \\ W & \cdots & W
\end{array} \right]. 
\] 

Then the evolution of the relative states is described by:
\begin{align*}
\dot x = A_n x + B_n u + \sqrt \epsilon F_n w_p  - G_n w_e, \quad x(0)=x_0,
\end{align*}
with observations
\[
y= C_n \, x + H_n v. 
\]
We will assume that $x_0$ is a Gaussian random vector with mean $\bar x_0$ and covariance matrix $\Sigma_0 \succ 0$. 

\subsection{Basic Example: Simple Integrators}

In an example used repetedly in the following, the pursuers are simply $n$ integrators and the evader follows a standard brownian motion:
\begin{align*}
\begin{cases} \dot x_e &= \; w_e \\
\dot x_{p,i} &= \; u_i + \sqrt{\epsilon} \; w_{p,i}, \quad i=1,\ldots,n,
\end{cases}
\end{align*}
i.e., we have the dynamics
\[
\dot x = u +\sqrt{\epsilon} \, w_p - \left[ \begin{array}{c} w_e \\ \vdots \\ w_e \end{array} \right].
\]
We will also take $C=H=I_d$, i.e.,
\[
y=c+v.
\]

%% file: risk.tex
\section{Risk-Sensitive Tracking} \label{section: solution}

Let $Q_n=I_n \otimes Q$ and $R_n=I_n \otimes R$. We formulate the tracking problem as an infinite horizon linear exponential quadratic Gaussian (LEQG) problem:
\begin{align*}
&\text{minimize } \\  J&=\lim_{T \to \infty} \frac{2}{\theta T} \ln \mathbb{E} \left\{  \exp \frac{\theta}{2} \int_0^T \frac{1}{n} \left( \sum_{i=1}^n (x_{p,i}-x_e)' Q (x_{p,i}-x_e) \right. \right. \\
& \left. \left.  + \sum_{i=1}^n u_{p,i}' \, R \, u_{p,i} \right) dt \right\} \\
&= \lim_{T \to \infty} \frac{2}{\theta T} \ln \mathbb{E} \left\{  \exp \frac{\theta}{2} \int_0^T \frac{1}{n} \left( x' Q_n x + u' R_n u \right) dt \right\}.
\end{align*}

The $\frac{1}{n}$ factor is due to the fact that we want to obtain a measure of performance per agent. If $\theta > 0$ we consider risk-averse agents, if $\theta < 0$ the agents are risk-seeking and in the limit of $\theta=0$, the agents are risk-neutral and we recover the standard LQG formulation. There is a very large litterature on the LEQG problem, as well as on more general risk-sensitive control problems, that we cannot review here. Suffices to say that the LEQG problem was introduced in the engineering litterature and solved in the full information case by Jacobson (who did not use the logarithmic transformation) \cite{Jacobson73_RSC}. Links with differential games are emphasized already in that paper, and the relationship to $H^\infty$ control is considered in \cite{Glover88, Basar95_book, whittle90_RSbook}. The risk-sensitive state estimator is not the Kalman filter, and robustness properties of this filter have been explored for example in \cite{Boel02_RSFiltering}. An important motivation for the early work of Jacobson in the full information case was to obtain a controller which is not independent of the noise in the dynamics, as arises in the LQG case because of the certainty equivalence principle. This motivation translates in our case into obtaining controllers for the pursuers which are not independent of the known characteristics of the random motion of the evader, certainly a desirable feature.

The LEQG output feedback solution over an infinite horizon is described for example in \cite{whittle90_RSbook, Pan96_model}. We start by considering the full state feedback problem or perfect measurement case, that is, $y=x$. Introduce the quantity 
\[
S_n(\theta)=n B_n R_n^{-1} B_n'-\theta \, (W_n+ \epsilon Z_n)
\]
and consider the generalized algebraic Riccati equation (GARE)
\begin{equation}	\label{eq: full state GARE}
A_n'X+X A_n-X \, S_n(\theta) \, X + \frac{1}{n} Q_n = 0.
\end{equation}
Now define the quantity
\begin{align}
\theta^*(n)=\sup \left\{ \theta \in \mathbb{R}: \text{ the GARE (\ref{eq: full state GARE}) admits a} \right. \nonumber \\
\left. \text{positive definite solution } X_{n,\theta} \right\}. \label{eq: critical value full info}
\end{align}

Assume that $(A,B)$ is controllable and $(A,Q)$ is observable. Then we have that $\theta^*(n)$ is positive, and for all $\theta < \theta^*(n)$, the LEQG problem with perfect state measurements admits an optimal state-feedback solution
\begin{equation} \label{eq: control}
u=-(n R_n^{-1} B_n'X_{n,\theta}) \, x,
\end{equation}
with the optimal cost being
\begin{equation}	\label{eq: full info cost}
J^*(\theta,n)=\Tr((W_n + \epsilon Z_n) \, X_{n,\theta}).
\end{equation}
Furthermore the feedback matrix $A_n-n B_n R_n^{-1} B_n' X_{n,\theta}$ is Hurwitz. In the following we will omit to write the dependence with respect to $\theta$ of the solution of the GARE, and write just $X_n$.

Now consider the output feedback solution for the problem with imperfect measurements described earlier. Recall that we assumed $V > 0$ and let $V_n=I_n \otimes V$. Introduce the quantity
\begin{align*}
T_n(\theta) = C'_n V_n^{-1} C_n - \theta \frac{Q_n}{n}, 
\end{align*}
and consider the two GAREs (\ref{eq: full state GARE}) and
\begin{equation}	\label{eq: output feedback GARE}
Y A_n'+A_n Y - Y \, T_n(\theta) \, Y + (W_n+ \epsilon Z_n)= 0.
\end{equation}
Define the quantity
\begin{align*}
& \theta_I^*(n)=\sup \{ \theta \in \mathbb{R}:  \text{ the GAREs (\ref{eq: full state GARE}) and (\ref{eq: output feedback GARE}) admit minimal } \\ 
&\text{ positive definite solutions } X_n \text{ and } Y_n, \text{ respectively, and } \\
& \text{ the matrix } I-\theta Y_n X_n \text{ has only positive eigenvalues} \}.
\end{align*}

For $\theta < \theta_I^*(n)$, introduce the filter
\begin{align*}
& \frac{d {\hat x}}{dt} = (A_n+\theta Y_n \frac{Q_n}{n}) \, \hat x  + B_n \, u  + Y_n C'_n V_n^{-1} (y - C_n \, \hat x ), \\
& \; \hat x_0=\bar x_0.
\end{align*}
Let $\tilde x=(I-\theta Y_n X_n)^{-1} \hat x$. Then one can compute that that $\tilde x$ is generated by the following differential equation:
\begin{align*}
& \frac{d \tilde x}{dt} = (A_n - S_n(\theta) X_n) \, \tilde x + (I-\theta Y_n X_n)^{-1} \, B \, \tilde u \\
& + (I-\theta Y_n X_n)^{-1} Y_n C'_n V_n^{-1} (y - C_n \tilde x), 
\end{align*}
where $\tilde u = u \, + \, n R_n^{-1} B_n' X_n \, \tilde x$.

Now suppose that the pairs $(A,B)$ and $(A_n, [\sqrt{\epsilon} F_n ,\; -G_n])$ are controllable. For the second condition, in practice we will assume $(A,F)$ to be controllable. Also, assume that the pairs $(A,C)$ and $(A,Q)$ are  observable, and $\Sigma_0 \leq Y_n$, then for all $\theta < \theta_I^*(n)$, the optimal controller is given by
\[
u^*=-n R_n^{-1} B'_n X_n (I-\theta Y_n X_n)^{-1} \, \hat x = -n R_n^{-1} B'_nX_n \, \tilde x.
\]
The optimal cost is
\begin{equation}	\label{eq: total cost}
J^*_I(\theta,n)=\text{Tr}\left( \, Y_n \frac{Q_n}{n}+Y_n C'_n V_n^{-1} C_n Y_n \, X_n \, (I-\theta Y_n X_n)^{-1}\, \right).
\end{equation}

\subsection{Application to the basic example}

\begin{figure} 
\centering
\subfigure[Risk Seeking Agents.] 
{
    \label{fig:sub:a}
    \includegraphics[width=4cm]{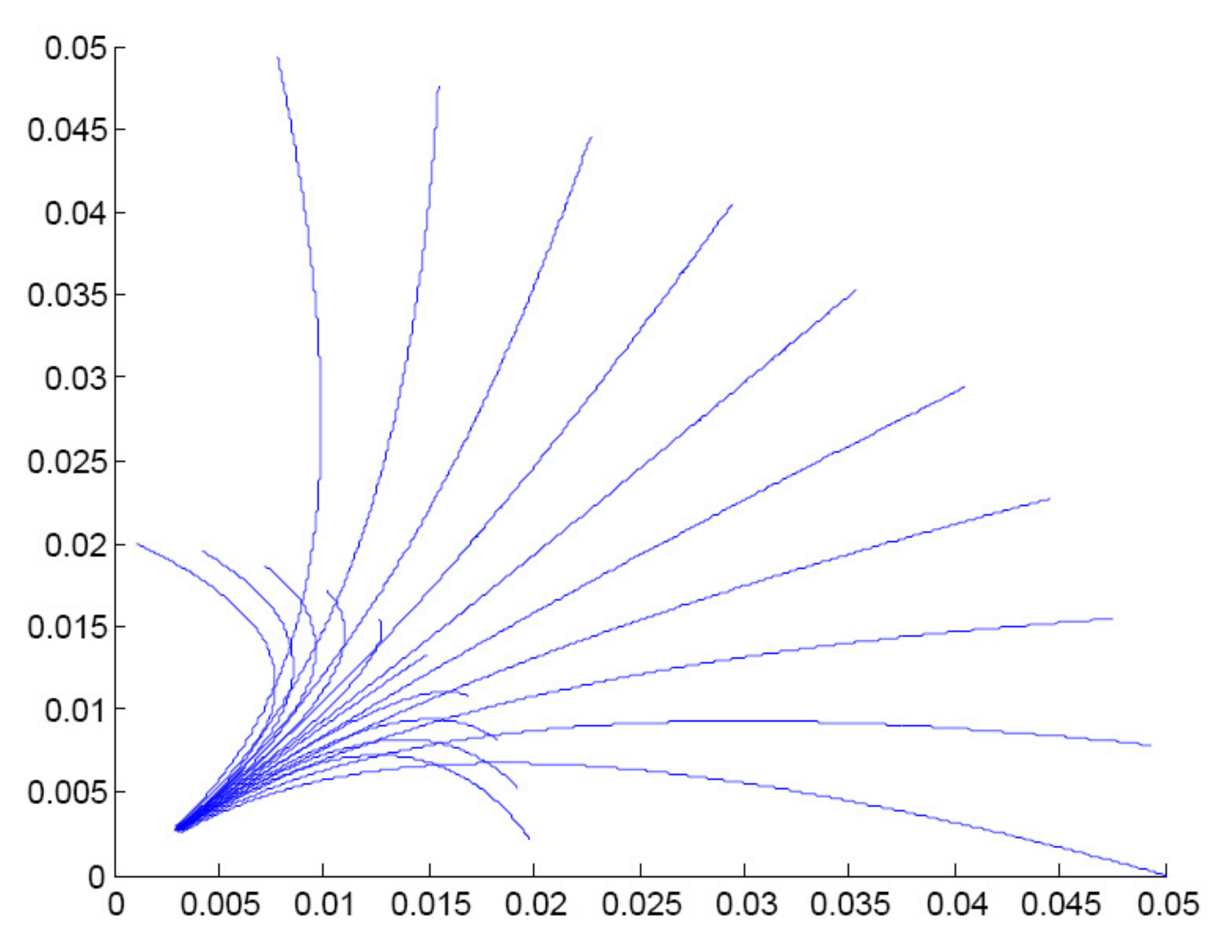}
}
\subfigure[Risk Neutral Agents.] 
{
    \label{fig:sub:b}
    \includegraphics[width=4cm]{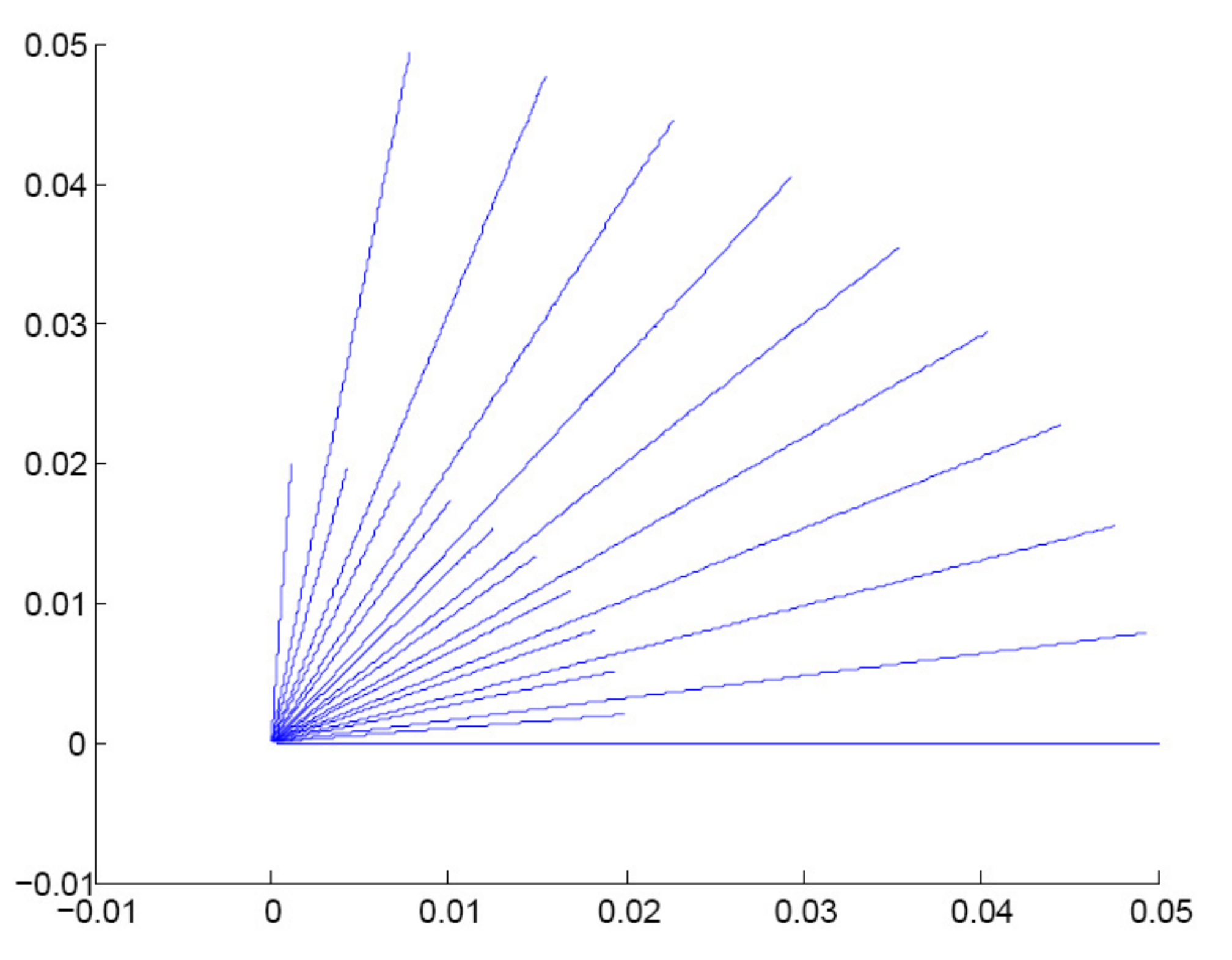}
}
\hspace{1cm}
\subfigure[Risk Averse Agents.] 
{
    \label{fig:sub:c}
    \includegraphics[width=4cm]{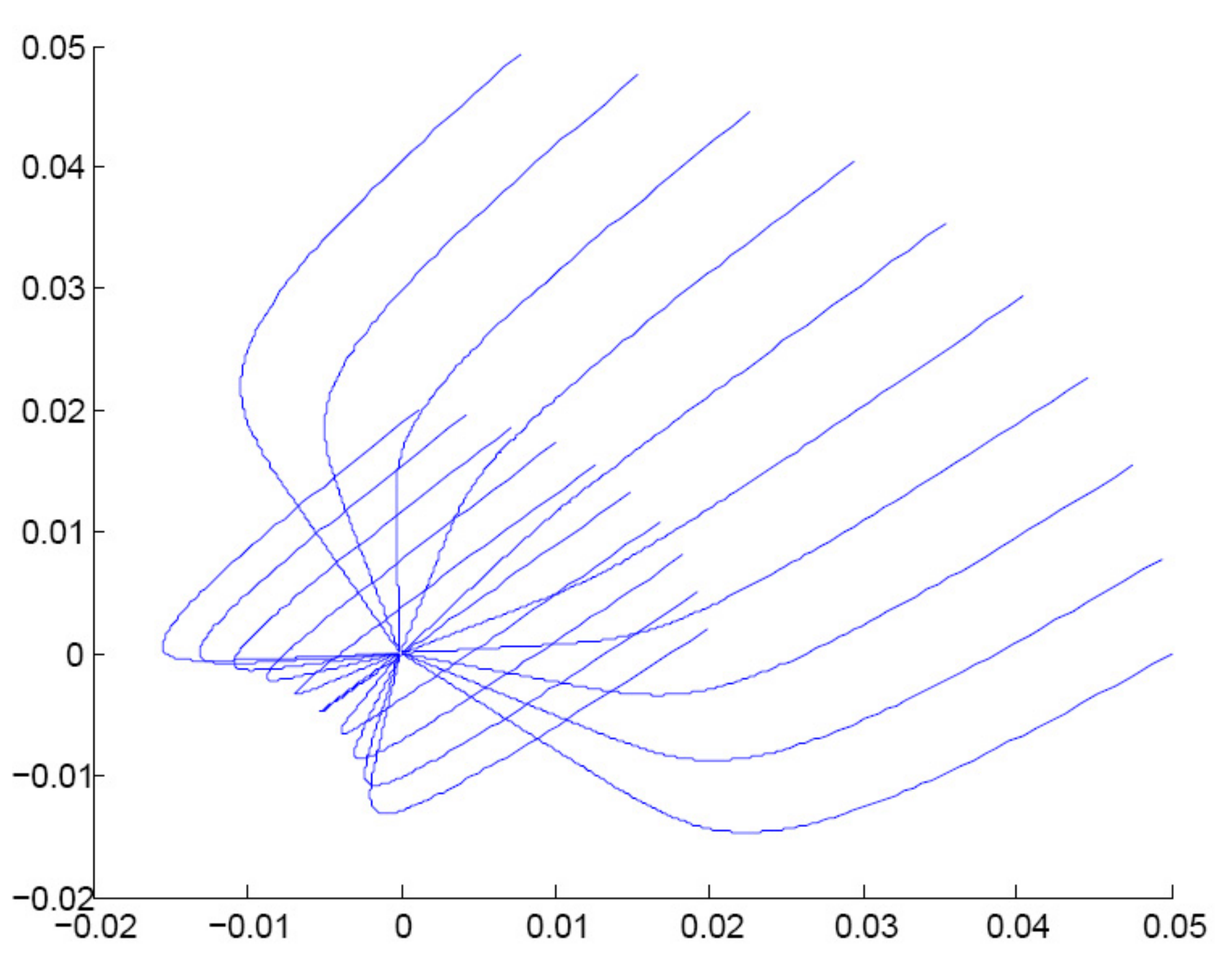}
}
\caption{Trajectories for Risk Sensitive Agents.}
\label{fig: trajectories} 
\end{figure}

It is interesting to visualize the trajectories of the agents in a simple case of the basic example. We consider a situation with perfect information, and $\epsilon=0$. The agents are tracking an evader according to the model of random motion described earlier, but in this simulation the evader in fact remains immobile at the origin. The trajectories are shown on Fig. \ref{fig: trajectories}, for risk-averse, risk-neutral and risk-seeking agents, and show different qualitative convergence behaviors.

The risk-neutral trajectories are essentially trivial in the sense that the controllers are totally decoupled and the control law for one agent depends only on its separation from the target. In the risk-seeking and risk-averse case however, the controllers become coupled, and each agent needs information about the position of all the other pursuers as weel as the target. In the risk averse case, two set of agents starting at different distances from the target try to track it from both sides. Intuitively, the overshoot is due to a pessimistic behavior which leads them to give more importance to the case where the target moves away from them.

%% file: performance.tex
\section{Performance Analysis} \label{section: performance}

Equations (\ref{eq: full info cost}) and (\ref{eq: total cost}) give the cost per agent when the group of pursuers consist of $n$ agents. In this section we study how this individual cost evolves as the number of agents increases.

\subsection{Perfect State Measurements}

\vspace{0.1cm}
\subsubsection{LQG solution}

The LQG solution is obtained for $\theta=0$ (risk neutral). In this case we can rewrite the control GARE (\ref{eq: full state GARE}) as:
\begin{align}
& A_n'X+X A_n-X \, n B_n R_n^{-1} B'_n \, X + \frac{1}{n} Q_n = 0 \label{eq: LQG control ARE}.
\end{align}

Due to the block diagonal form of the matrices, the equation decouples by agents. With our controllability and observability assumptions, (\ref{eq: LQG control ARE}) has a unique positive definite solution (\cite{zhou96_book}, corollary 13.8). In particular, if $X_1$ is the positive definite solution for one agent, i.e.,
\[
A' X_1+X_1 A - X  B R^{-1} B' \, X + Q = 0
\]
then by unicity the solution for $n$ agents is verified to be
\[
X_n=\frac{1}{n}I_n \otimes X_1.
\]
The controller for the group, $u=- (I_n \otimes R^{-1} B' X_1) \, x$ is decoupled as well, that is, each agent considers only the state difference between himself and the target, but not the states of the other agents. 

The cost per agent in the perfect information case is then
\begin{align*}
J^*(0,n) &= \Tr((W_n + \epsilon \, Z_n) \, X_n) = \frac{1}{n} \Tr(E_n \otimes W X_1 + I_n \otimes \epsilon \, Z X_1) \\
&= \Tr((W + \epsilon \, Z) X_1) = J^*(0,1),
\end{align*}
that is, the cost per agent is independent of the number of agents.

\vspace{0.1cm}
\subsubsection{Risk-Sensitive Solution with no dynamics noise}

When $\theta \neq 0$, some coupling between the controllers of the agents is introduced through the matrix $W_n$. This is a desirable feature as intuitively we would like the agents to take advantage of the fact that they can cooperate. We will take $\epsilon=0$ here, which still leads to a nonsingular problem in the perfect information case. 

Consider first the case of one agent. We can always write the solution as a pertubation of the LQG solution, i.e., as
\[
X_1=\tilde X_1+\hat X_1,
\]
where $X_1$ is the solution to (\ref{eq: full state GARE}) for $n=1$ and $\tilde X_1$ is the solution to (\ref{eq: LQG control ARE}) for $n=1$. 

\vspace{0.1cm}

\begin{proposition}
For $n$ agents and $\epsilon=0$, the LEQG problem with perfect state measurements admits an optimal state-feedback solution solution (\ref{eq: control}) with the solution to (\ref{eq: full state GARE}) given by
\begin{equation}	\label{eq: solution LQG full state}
X_n= \frac{1}{n} I_n \otimes \tilde X_1 + \frac{1}{n^2} E_n \otimes \hat X_1.
\end{equation}
\end{proposition}

\vspace{0.3cm}

\begin{proof}
By definition we have that $\tilde X_1$ and $\hat X_1$ verify:
\begin{align}
& A' \tilde X_1 + \tilde X_1 A - \tilde X_1 B R^{-1} B \tilde X_1 + Q = 0 \label{eq: interm1} \\
& A' (\tilde X_1 + \hat X_1) + (\tilde X_1 + \hat X_1) A \nonumber \\
& - (\tilde X_1 + \hat X_1) (B R^{-1} B' - \theta W) (\tilde X_1 + \hat X_1) + Q = 0. \label{eq: interm2}
\end{align}

We try this candidate solution in the GARE (\ref{eq: full state GARE}) for $\epsilon=0$. Then by a straightforward calculation, using the identity $E_n^2 = n E_n$, we obtain:
\begin{align*}
& \frac{1}{n} I_n \otimes (A' \tilde X_1 + \tilde X_1 A - \tilde X_1 B R^{-1} B \tilde X_1 + Q) + \\
& \frac{1}{n^2} E_n \otimes \left[ A' (\tilde X_1 + \hat X_1) + (\tilde X_1 + \hat X_1) A \right. \nonumber \\
& - (\tilde X_1 + \hat X_1) (B R^{-1} B' - \theta W) (\tilde X_1 + \hat X_1) \\
& \left. - (A' \tilde X_1 + \tilde X_1 A - \tilde X_1 B R^{-1} B) \right] \\
& = 0.
\end{align*}
Indeed, in the first term we have (\ref{eq: interm1}) and in the second term we have the difference between (\ref{eq: interm2}) and (\ref{eq: interm1}).

Moreover, this solution $X_n$ is positive definite and stabilizing. This follows from the fact that the eigenvalues of
\[
I_n \otimes M + \frac{E_n}{n} \otimes N,
\]
for any diagonalizable matrices $M$ and $(M+N)$, are the eigenvalues of $M$ with multiplicities $(n-1)$ as well as those of $(M+N)$. A set of corresponding eigenvectors are $y_i \otimes w_j$, where the $y_i$ are $(n-1)$ vectors spanning the kernel of $E_n$ (which is symmetric) and the $w_j$ are eigenvectors of $M$, and $\mathbf{1}_n \otimes v_i$, for $v_i$ eigenvector of $(M+N)$. It is then easy to see that the eigenvalues of $X_n$ are thoses of $X_1$ and $\tilde X_1$, and those of $A_n-n B_n R_n^{-1} B_n' X_n$ are the eigenvalues of $A-B R^{-1} B' X_1$ and $A-B R^{-1} B' \tilde X_1$. But $X_1$ and $\tilde X_1$ are stabilizing and positive definite.
\end{proof}

\vspace{0.1cm}
\begin{remark}
In particular, we see that the critical value $\theta^*(n)$ is independent of $n$ in the perfect information case when $\epsilon=0$. 
\end{remark}
\vspace{0.1cm}

We can then compute the cost of this solution to the risk-sensitive tracking problem, for $\epsilon=0$. We get
\begin{align*}
J^*(\theta,n) &= \Tr(W_n  \, X_n) = \Tr(\frac{1}{n} E_n \otimes W \tilde X_1 + \frac{1}{n^2} E_n^2 \otimes W \hat X_1) \\
&= \Tr(W (\tilde X_1 + \hat X_1)) = J^*(\theta,1).
\end{align*}
Hence, as in the LQG case, the cost per agent in the case $\epsilon=0$ is independent of the number of agents.

\vspace{0.1cm}
 
\subsubsection{Risk-Sensitive Solution with noisy dynamics}

The risk-sensitive solution in the perfect measurement case seems to be most relevant when noise is present in the dynamics of the agents. In this case, the cost-per-agent is not independent of the number of agent any more, and moreover, the critical value of the parameter $\theta^*(n)$ increases with $n$. Another way to say this is that for a fixed value of risk aversion $\theta^*$, there is a minimum number of agents that are necessary to obtain a finite cost.

\begin{figure}[thpb]
\centering
\includegraphics[height=9cm]{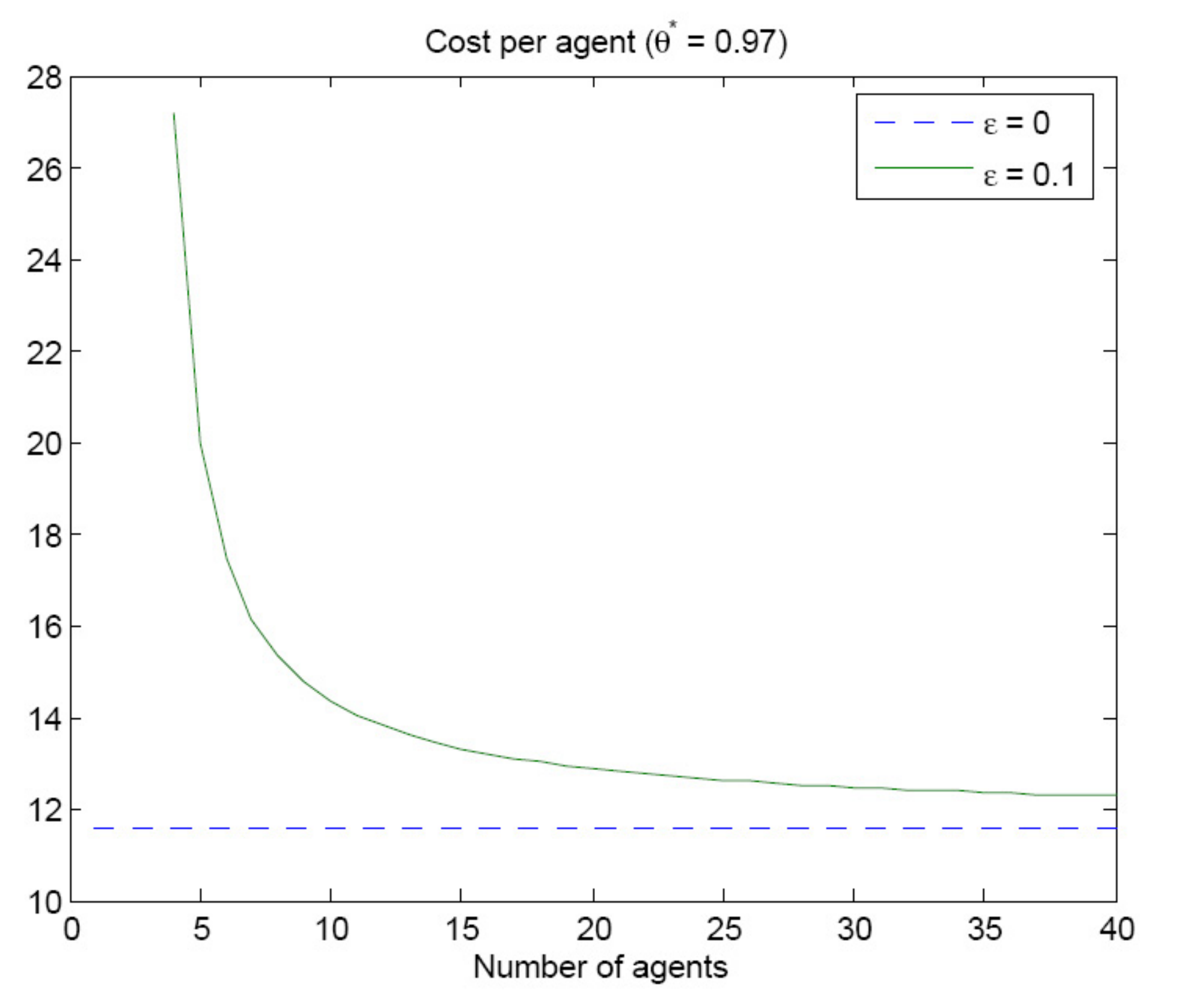} 
\caption{Cost per agent in the basic example with noisy dynamics, under perfect measurements. The cost is infinite for less than 4 agents.}
\label{fig: performance perfect info}
\end{figure}

We only provide a numerical illustration of this fact. Fig. \ref{fig: performance perfect info} shows the cost per agent in the basic example with $\epsilon=0.1$. We set $\theta^*=0.97$. In this case, the cost per agent is finite only if $4$ agents or more are used. The lower bound is the cost per agent for $\epsilon=0$.

\input{imperfect}

%% file: imperfect.tex
\subsection{Imperfect State Measurements}

In the case of imperfect measurements, there is one obvious advantage of using more agents. As the number of agents taking measurements with independent noise sources increases, a better state estimate can be constructed. In this section, we consider the imperfect measurement tracking problem, with the simplifying assumption $A=0$. Note that in this case, the matrix $(A_n,\sqrt{\epsilon} F_n - G_n)$ is not stabilizable when $\epsilon=0$, resulting in an only marginally stable filter. Hence we have to keep a small noise term in the calculation, but in the analysis we will consider the dominant part as $\epsilon \to 0$, for which a clear answer is available.

\vspace{0.1cm}
\subsubsection{LQG solution}

For $\theta=0$ (risk neutral), the filter algebraic Riccati equation (\ref{eq: output feedback GARE}) for $A=0$ becomes:
\begin{align}
& - Y \, C'_n V_n^{-1} C_n \, Y + (W_n + \epsilon Z_n) = 0. \label{eq: LQG filter ARE}
\end{align}

Let $Y_1$ be the positive definite solution (unique in the LQG case) of the Riccati equation (\ref{eq: LQG filter ARE}) for one agent, i.e.,
\[
- Y_1 \, C' V^{-1} C \, Y_1 + (W + \epsilon Z)= 0.
\]

Then we have
\begin{align*}
& (\frac{E_n}{\sqrt{n}} \otimes Y_1) (I_n \otimes C V^{-1} C) (\frac{E_n}{\sqrt{n}} \otimes Y_1) = (\frac{E_n^2}{n}) \otimes (W+\epsilon Z) \\
& = E_n \otimes W + \epsilon{E_n \otimes Z},
\end{align*}
and so as $\epsilon \to 0$, the solution to the $n$ agent problem, for which the right hand side in the previous equation is $E_n \otimes W + \epsilon{I_n \otimes Z}$, approaches $Y_n= \frac{E_n}{\sqrt{n}} \otimes Y_1$. 

The total cost is
\[
J^*_I(0,n)=\Tr(Y_n \frac{Q_n}{n}+(W_n+\epsilon Z_n) X_n).
\]
We have already seen that if $\epsilon=0$, the second term of this expression becomes independent of $n$. Hence, as $\epsilon \to 0$, the cost approaches
\begin{align*}
J^*_I(0,n) & \approx \frac{1}{n} \text{Tr}\left\{ (\frac{E_n}{\sqrt{n}} \otimes Y_1) (I_n \otimes Q) \right\} + \Tr(W X_1) \\
& \approx \frac{1}{\sqrt{n}} \text{Tr} \left\{ Y_1 Q \right\} + \Tr(W X_1).
\end{align*}
In conclusion, as the number of agents increases, the tracking performance per agent converges to the control performance for one agent at rate $1/\sqrt{n}$, due to a better estimation performance only. This is intuitively expected from our understanding of the asymptotic normality of maximum likelihood estimators. However, if we consider the more general diffusion process with $A \neq 0$ and $\epsilon \neq 0$, it is not clear if this asymptotic rate of convergence still holds.

\vspace{0.1cm}
\subsubsection{Risk-Sensitive Solution}

Again let us consider the case $A=0$ and $\epsilon \to 0$. In the limit $\epsilon=0$, the solution to the control GARE is given by (\ref{eq: solution LQG full state}) and this equation does not influence the critical value of the parameter $\theta_I^*$. The filter GARE is:
\begin{align}
& Y (I_n \otimes (C V^{-1} C - \frac{\theta}{n} Q)) Y = E_n \otimes W + \epsilon I_n \otimes Z \label{eq: RS estimation ARE II}. 
\end{align}

\begin{figure}[thpb]
\centering
\includegraphics[height=9cm]{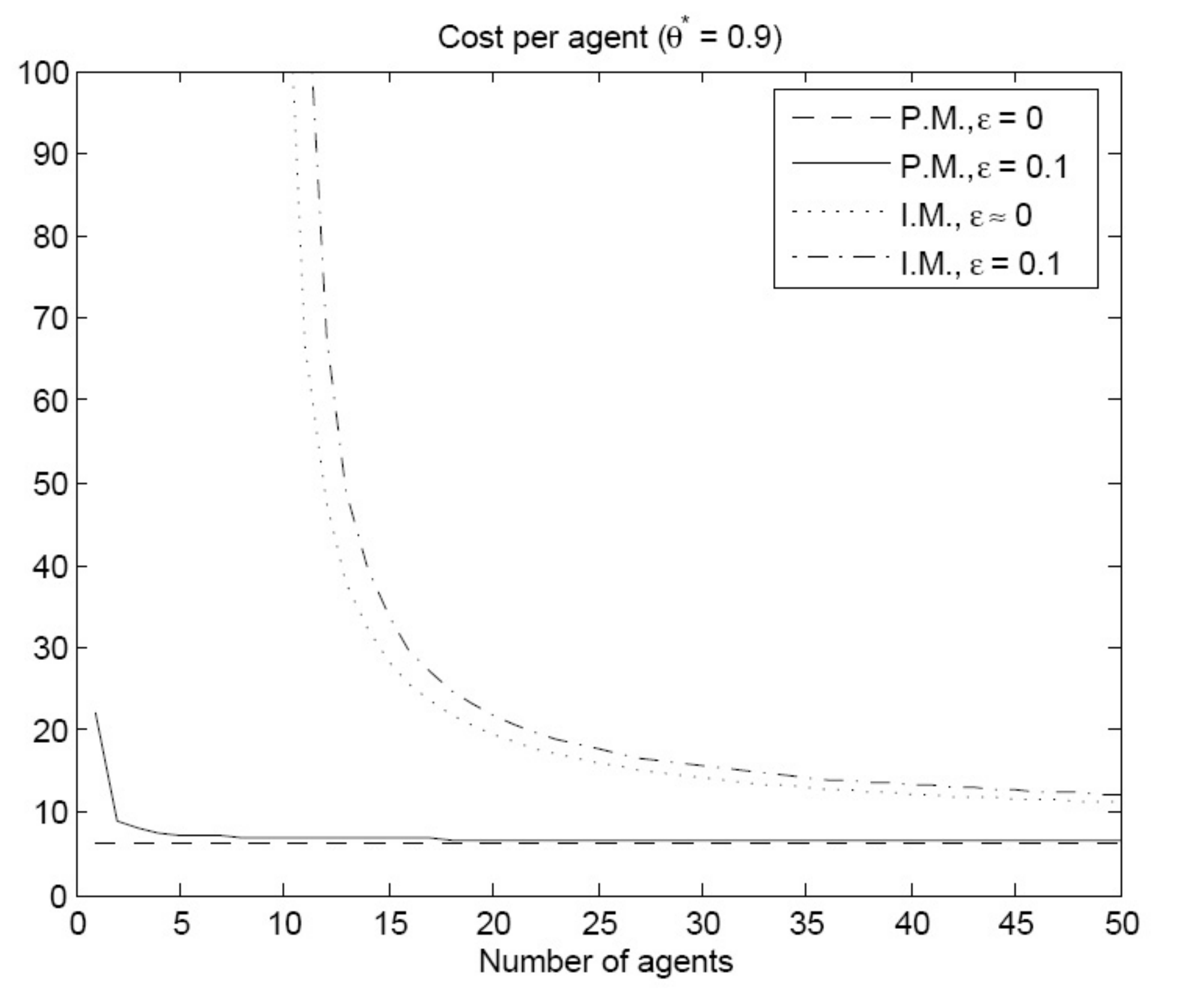} 
\caption{Cost per agent in the basic example with noisy dynamics and observations (P.M.=perfect measurements, I.M.=imperfect measurements. A=0). In the case $\epsilon=0.1$ and noisy observations, the cost is infinite for less than 10 agents.}
\label{fig: performance imperfect info}
\end{figure}

Hence as in the LQG problem, we see that as $\epsilon \to 0$, the solution to this equation approaches
\begin{equation}	\label{eq: solution filter}
Y_n = \frac{E_n}{\sqrt{n}} \otimes \tilde Y_{1,n},
\end{equation}
but an essential difference is that now, $\tilde Y_{1,n}$ is the solution to the equation:
\begin{equation} \label{eq: solution filter small}
Y (C V^{-1} C - \frac{\theta}{n} Q) Y = W+\epsilon Z,
\end{equation}
which is \emph{not} the single agent equation, but \emph{depends on $n$} through the parameter $\theta/n$. Fig. \ref{fig: performance imperfect info} provides some experimental results on the performance of the multi-agent system in various cases considered above.

\begin{figure}[thpb]
\centering
\includegraphics[height=9cm]{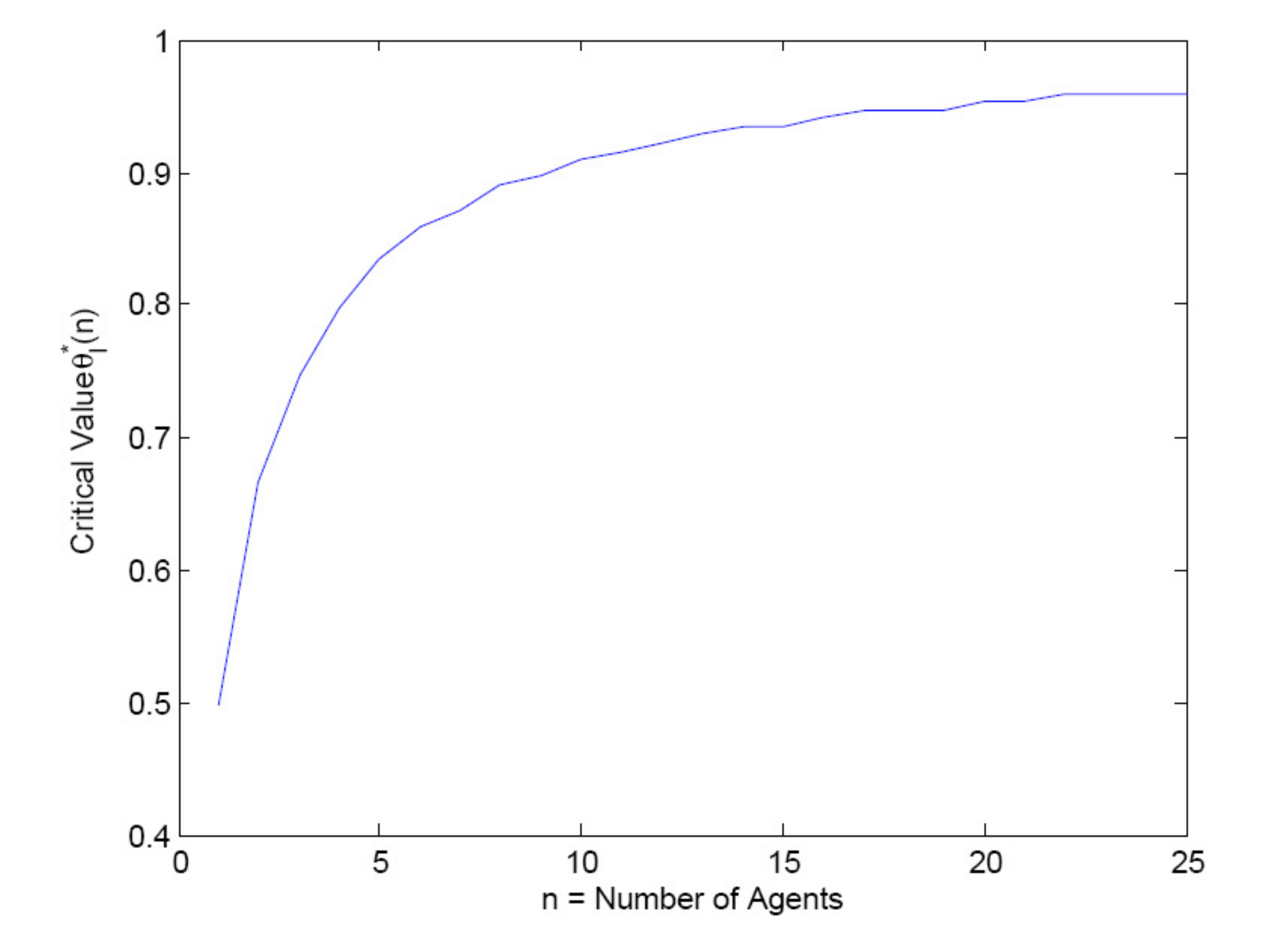} 
\caption{Evolution of $\theta_I^*(n)$ with the number of agents in the basic example ($\epsilon \approx 0$, A=0).}
\label{fig: theta star}
\end{figure}

If the constraint that $I-\theta Y_n X_n$ must have only positive eigenvalues were not present, from (\ref{eq: solution filter small}) $\theta_I^*(n)$ would increase linearly with $n$. Now it is easy to check from (\ref{eq: solution LQG full state}) and (\ref{eq: solution filter}) that as $\epsilon \to 0$, $Y_n X_n$ approaches
\[
Y_n X_n \approx \frac{E_n}{n^{3/2}} \otimes \tilde Y_{1,n} X_1,
\]
$E_n/n$ has eigenvalues $1$ and $0$, so the condition as $\epsilon \to 0$ becomes 
\[
\rho(\theta \tilde Y_{1,n,\theta} X_{1,\theta}) < \sqrt{n},
\]  
where $\rho$ denotes the spectral radius. We have indicated the dependence of $Y_{1,n}$ and $X_1$ on $\theta$ in this condition. Fig. $\ref{fig: theta star}$ shows the evolution of the critical value $\theta_I^*(n)$ with the number of agents, computed for the basic example. We obtained similar experimental results for $\epsilon \neq 0$. However if $A \neq 0$ it seems that we can obtain different behaviors.

Note that in the various cases considered above, when an analytical solution could be obtained it always implied that for solving the n agent problem, we only need to solve a Riccati equation of the same size as for one agent. Hence significant computational reduction can be achieved by taking advantage of the symmetries present when considering homogeneous agents.

